\documentclass[11pt]{amsart}

\usepackage{latexsym}
\usepackage{psfrag}
\usepackage{amsmath}
\usepackage{amssymb}
\usepackage{epsfig}
\usepackage{amsfonts}

\def\wh{\widehat}

\def\zz{\mathbb Z}

\def\sm{\smallsetminus}
\def\pr{\prime}

\def\la{\lambda}
\def\ga{\gamma}
\def\si{\sigma}

\def\al{\alpha}
\def\be{\beta}

\def\vp{\varphi}

\def\ch{\mathcal H}

\def\cp{\mathcal P}
\def\cq{\mathcal Q}

\def\<{\langle}
\def\>{\rangle}

\newtheorem{thm}{Theorem}[section]
\newtheorem{lemma}[thm]{Lemma}

\newtheorem{Remark}[thm]{Remark}

\begin{document}

\title[A combinatorial proof of the Rogers-Ramanujan and Schur identities]%
    {A combinatorial proof of the Rogers-Ramanujan and Schur identities}
\author[Cilanne Boulet]{Cilanne Boulet$^\ast$}
\author[Igor~Pak]{Igor~Pak$^\ast$}
\date{September 12, 2005}


\thanks{${\hspace{-.95ex}}^\ast$Department of Mathematics, MIT, Cambridge, MA 02139, 
\ Email: \texttt{\{cilanne,pak\}@math.mit.edu}}

\begin{abstract}
\noindent
We give a combinatorial proof of the first Rogers-Ramanujan
identity by using two symmetries of a new generalization
of Dyson's rank.  These symmetries are established
by direct bijections.
\end{abstract}

\maketitle

\section*{Introduction} \label{sec:intro}

The Roger-Ramanujan identities are perhaps the most mysterious and
celebrated results in partition theory.  They have a remarkable
tenacity to appear in areas as distinct as enumerative combinatorics,
number theory, representation theory, group theory,
statistical physics, probability and complex analysis~\cite{A1,A4}.
The identities were discovered independently by Rogers, Schur, and
Ramanujan (in this order), but were named and publicized by
Hardy~\cite{H}. Since then, the identities have been greatly
romanticized and have achieved nearly royal status in the field.
By now there are dozens
of proofs known, of various degree of difficulty and depth.
Still, it seems that Hardy's famous comment remains valid:
{\em ``None of the proofs of} \, [the Rogers-Ramanujan identities]
{\em can be called ``simple'' and ``straightforward''} [...];
{\em and  no doubt it would be unreasonable to expect a
really easy proof''}~\cite{H}.

In this paper we propose a new combinatorial proof
of the first Rogers-Ramanujan identity with a minimum
amount of algebraic manipulation.
Almost completely bijective, our proof would not satisfy
Hardy as it is neither ``simple" nor ``straightforward''.
On the other hand, the heart of the proof is the analysis
of two bijections and their properties, each of them
elementary and approachable.  In fact, our proof gives
new generating function formulas (see~$(\maltese)$ in Section~1)
and is amenable to advanced generalizations which will appear
elsewhere (see~\cite{Bo}).

We should mention that on the one hand, our proof is heavily
influenced by the works of Bressoud and Zeilberger~\cite{B,BZ1,BZ2,BZ3},
and on the other hand by Dyson's papers~\cite{D1,D2}, which were
further extended by Berkovich and Garvan~\cite{BG}
(see also~\cite{G,P2}).
In fact, the basic idea to use a generalization of Dyson's
rank was explicit in~\cite{BG,G}.  We postpone historical and
other comments until Section~\ref{sec:final}.

Let us say a few words about the structure of the paper.
We split the proof of the first Rogers-Ramanujan identity into two
virtually independent parts.  In the first, {\em the algebraic part},
we use the Jacobi triple product identity
to derive the identity
from two \emph{symmetry equations}.  The latter  are proved in
{\em the combinatorial part} by direct bijections.
Our presentation is elementary and completely self-contained,
except for the use of the classical Jacobi triple product
identity.  We conclude with the final remarks section.

\bigskip


\section{The algebraic part} \label{sec:alg}

We consider  the first \emph{Rogers-Ramanujan identity}:
$$(\blacklozenge) \, \qquad
1+\sum_{k=1}^\infty \frac{t^{k^2}}{(1-t)(1-t^2)\cdots (1-t^k)} \ = \
\prod_{i=0}^\infty \frac{1}{(1-t^{5i+1})(1-t^{5i+4})} \, .$$

Our first step is standard.  Recall the \emph{Jacobi triple
product identity} (see e.g.~\cite{A1}):
$$\sum_{k=-\infty}^\infty z^k q^{\frac{k(k+1)}{2}}  \,
= \, \prod_{i=1}^\infty (1+zq^i) \,
\prod_{j=0}^\infty (1+z^{-1}q^{j}) \,
\prod_{i=1}^\infty (1-q^i).
$$
Set $q \gets t^5$, $z \gets (-t^{-2})$ and rewrite the right hand side of
$(\blacklozenge)$ as follows:
$$
\prod_{r=0}^\infty \frac{1}{(1-t^{5r+1})(1-t^{5r+4})} \ = \
\sum_{m=-\infty}^\infty \, (-1)^m \,t^{\frac{m(5m-1)}{2}} \
\prod_{i=1}^\infty \frac{1}{(1-t^i)}\,.
$$
This gives us \emph{Schur's identity}, which is equivalent to~$(\blacklozenge)$\,:
$$(\lozenge)  \quad
\left(1+\sum_{k=1}^\infty \frac{t^{k^2}}{(1-t)(1-t^2)\cdots (1-t^k)}\right)
\ = \  \prod_{i=1}^\infty \frac{1}{(1-t^i)} \
\sum_{m=-\infty}^\infty \, (-1)^m \, t^{\frac{m(5m-1)}{2}}
\,.
$$

\smallskip

To prove Schur's identity we need several combinatorial definitions.
Denote by $\cp_n$ the set of all partitions~$\la$ of~$n$, and
let~$\cp = \cup_n \cp_n$, $p(n) = |\cp_n|$.
Denote by~$\ell(\la)$ and~$e(\la)$
the number of parts and the smallest part of the partition,
respectively.  By definition,~$e(\la) = \la_{\ell(\la)}$.  We say
that $\la$ is a {\em Rogers-Ramanujan}  {\em partition} if
$e(\la) \ge \ell(\la)$.  Denote by $\cq_n$ the set of
Rogers-Ramanujan partitions, and let $\cq = \cup_n \cq_n$,
$q(n) = |\cq_n|$.  Recall that
$$P(t) \, := \, 1 \, + \, \sum_{n=1}^\infty \, p(n)\,t^n \, = \,
\prod_{i=1}^n \, \frac{1}{1-t^i}\,,
$$
and
$$Q(t) \, := \, 1 \, + \, \sum_{n=1}^\infty \, q(n)\,t^n \, = \,
1 \, + \,
\sum_{k=1}^\infty \frac{t^{k^2}}{(1-t)(1-t^2)\cdots (1-t^k)}\,.
$$

\smallskip

We consider a statistic on~$\cp \sm \cq$, the set of non-Rogers-Ramanujan partitions,
which we call the~{\em $(2,0)$-rank} of a partition, and denote by $r_{2,0}(\la)$,
for $\la \in \cp\sm \cq$.  Similarly, for~$m \ge 1$ we
consider a statistic on~$\cp$ which we call the~{\em $(2,m)$-rank}
of a partition, and denote by $r_{2,m}(\la)$, for $\la \in \cp$.
We formally define and study these statistics in the next section.
Denote by $h(n,m,r)$ the number of partitions~$\la$ of~$n$
with $r_{2,m}(\la) = r$. Similarly, let $h(n,m, \le r)$
and~$h(n,m, \ge r)$ be the number of partitions with the
$(2,m)$-rank~$\le r$ and~$\ge r$, respectively.
The following is apparent from the definitions:
$$(\divideontimes) \qquad \quad
\aligned
h(n,m,\leq r) \, + \, h(n,m,\geq r+1) & \, = \, p(n),
\ \,  \text{for} \ \, m>0, \ \text{and} \\
h(n,0,\leq r) \, + \, h(n,0,\geq r+1) & \, = \, p(n) - q(n),
\endaligned
$$
for all $r \in \zz$ and~$n\ge 1$.  The following two equations
are the main ingredients of the proof. We
have:

\smallskip

\textbf{(first symmetry)} \hskip2.36cm
$h(n,0,r) \, = \, h(n,0,-r)$,  and

\smallskip

\textbf{(second symmetry)}  \hskip1.cm
$h(n,m,\leq -r) \, = \, h(n-r-2m-2,m+2,\geq -r)$.

\smallskip
\noindent
The first symmetry holds for any $r$  and the second symmetry holds for~$m, r > 0$ and for~$m=0$ and~$r \geq 0$.

\noindent
Both symmetry equations will be proved in the next section.
For now, let us continue to prove Schur's identity.
For every $j \ge 0$ let
$$\aligned
a_j \, & = \, h\left(n-jr - 2jm - j(5j-1)/2, m+2j, \le -r -j\right), \ \text{ and}\\
b_j \, & = \, h\left(n-jr - 2jm - j(5j-1)/2, m+2j, \ge -r -j+1\right).
\endaligned
$$
The equation $(\divideontimes)$ gives us
$a_j + b_j = p(n -jr -2jm - j(5j-1)/2)$, for all
$r, j >0$.  The second symmetry equation gives us
$a_j = b_{j+1}$. Applying these multiple times
we get:
$$ \aligned
& h(n,m,\le -r) \, = \, a_0 \,  = \, b_1 \\
& \qquad \, =  \, b_1 + (a_1 - b_2) - (a_2 - b_3) + (a_3 - b_4) - \ldots \\
& \qquad \, =  \, (b_1 + a_1) - (b_2 + a_2) + (b_3 + a_3)  - (b_4 + a_4) + \ldots \\
& \qquad \, = \, p(n-r-2m-2) -p(n-2r-4m-9) +p(n-3r-6m-21) - \ldots  \\
& \qquad \, = \, \sum_{j=1}^\infty (-1)^{j-1} p(n -jr -2jm - j(5j-1)/2) \, .
\endaligned
$$
In terms of the generating functions
$$
H_{m,\le -r} (t) \, := \, \sum_{n=1}^\infty \, h(n,m,\le -r) \, t^n \,,
$$
this gives (for $m,r > 0$ and for $m=0$ and $r \geq 0$)
$$(\maltese) \qquad
H_{m,\le -r} (t) \, = \,
\prod_{n=1}^\infty \frac{1}{(1-t^n)} \ \sum_{j=1}^\infty (-1)^{j-1}
\, t^{jr+2jm+j(5j-1)/2}\,.
$$
In particular, we have:
$$\aligned
H_{0,\le 0} (t) \, & = \, \prod_{n=1}^\infty \frac{1}{(1-t^n)} \
\sum_{j=1}^\infty (-1)^{j-1} \, t^{\frac{j(5j-1)}2}\,, \\
H_{0,\le -1} (t) \, & = \, \prod_{n=1}^\infty \frac{1}{(1-t^n)} \
\sum_{j=1}^\infty (-1)^{j-1} \, t^{\frac{j(5j+1)}2}\,.
\endaligned
$$
From the first symmetry equation and~$(\divideontimes)$ we have:
$$H_{0,\le 0} (t) + H_{0,\le -1}(t) \, = \,
H_{0,\le 0} (t) + H_{0,\ge 1}(t) \, = \, P(t) - Q(t).$$
We conclude:
$$\aligned &
\prod_{n=1}^\infty \frac{1}{(1-t^n)}
\left(\sum_{j=1}^\infty \, (-1)^{j-1} \, t^\frac{j(5j-1)}{2} +
\sum_{j=1}^\infty \, (-1)^{j-1} \, t^\frac{j(5j+1)}{2}\right) \\
& \qquad = \,
\prod_{n=1}^\infty \frac{1}{(1-t^n)} \ -  \left(1
+ \sum_{k=1}^\infty \frac{t^{k^2}}{(1-t)(1-t^2) \dots (1-t^k)}\right),
\endaligned
$$
which implies~$(\lozenge)$ and completes the
proof of~$(\blacklozenge)$.

\bigskip


\section{The combinatorial part} \label{sec:comb}

\subsection{Definitions}\label{sec-def} 
\begin{figure}[hbt]
\begin{center}
\psfrag{la}{$\la$}
\psfrag{lp}{$\la^\pr$}
\epsfig{file=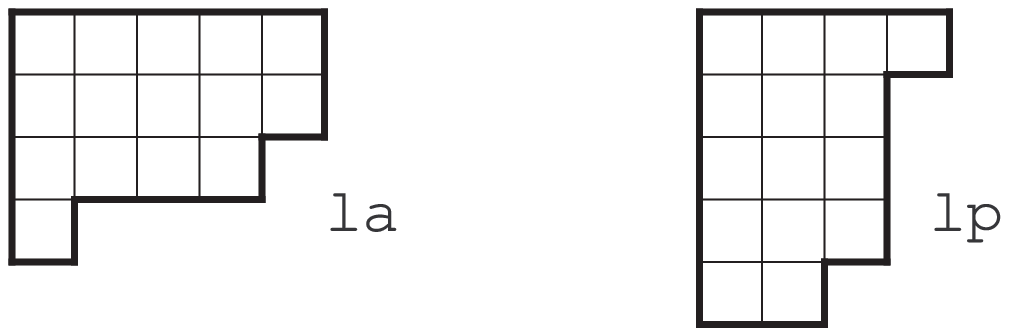,height=2.4cm}
\end{center}
\caption{Partition $\la = (5,5,4,1)$ and conjugate partition
$\la^\pr = (4,3,3,3,2)$. }
\label{fig:part}
\end{figure}

\
Let $\la = (\la_1,  \ldots,\la_{\ell(\la)})$, $\la_1 \ge \ldots \ge \la_{\ell(\la)} > 0$,
be an integer partition of~$n = \la_1 + \ldots + \la_{\ell(\la)}$.
We will say that $\la_j = 0$ for $j > \ell(\la)$.
We graphically represent the partition~$\la$
by a Young diagram~$[\la]$ as in Figure~\ref{fig:part}.
Denote by $\la^\pr$ the {\em conjugate partition} of~$\la$
obtained by reflection upon main diagonal
(see Figure~\ref{fig:part}).

\begin{figure}[hbt]
\begin{center}
\psfrag{la}{$\la$}
\psfrag{lp}{$\la^\pr$}
\psfrag{al}{$\al$}
\psfrag{be}{$\be$}
\psfrag{ga}{$\ga$}
\psfrag{a2}{$\ \al_2$}
\psfrag{b1}{$\, \be_1$}
\psfrag{g1}{$\ \ga_1^\pr$}
\epsfig{file=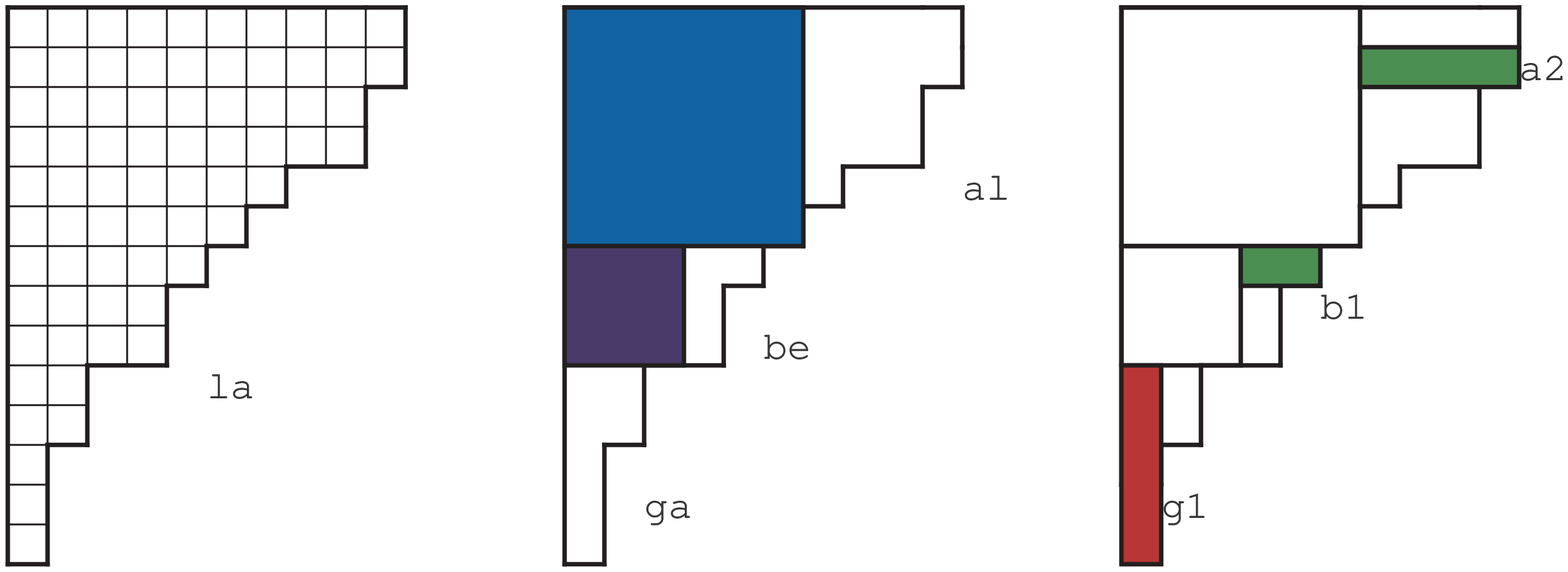,height=4.8cm}
\end{center}
\caption{Partition $\la = (10,10,9,9,7,6,5,4,4,2,2,1,1,1)$,
the first Durfee square of height~$s_0(\la) = 6$, and the
second Durfee square of height~$t_0(\la) = 3$.  Here the
remaining partitions are $\al = (4,4,3,3,1)$, $\be = (2,1,1)$,
and $\ga = (2,2,1,1,1)$. In this case, the $(2,0)$-rank
is $r_{2,0}(\la) = \be_1 + \al_2 - \ga_1^\pr = 2+4-5= 1$. }
\label{fig:durfee}
\end{figure}

\begin{figure}[hbt]
\begin{center}
\psfrag{la}{$\la$}
\psfrag{al}{$\al$}
\psfrag{be}{$\be$}
\psfrag{a2}{$\ \al_2$}
\psfrag{b1}{$\, \be_1$}
\epsfig{file=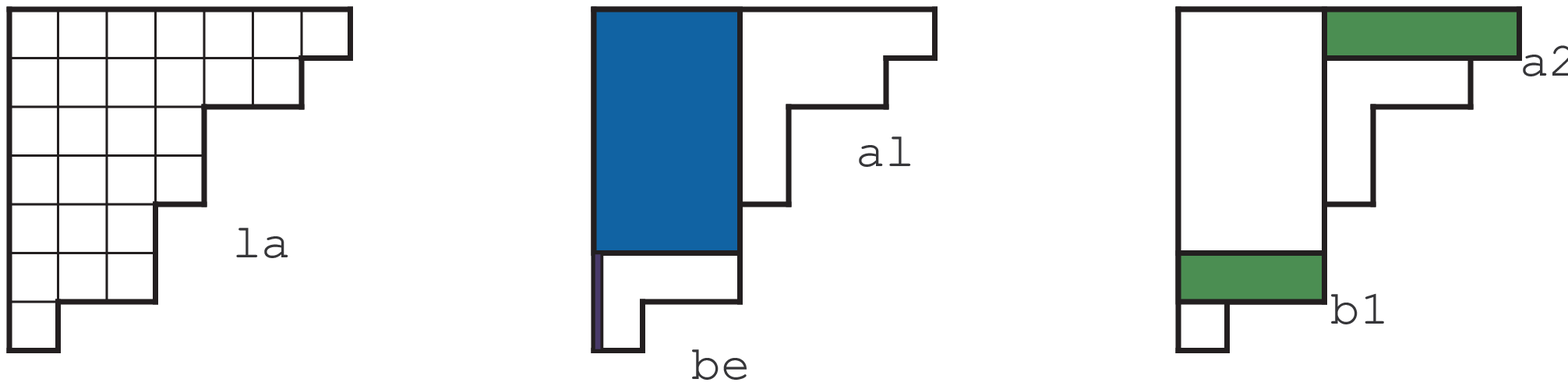,height=3.2cm}
\end{center}
\caption{Partition $\la = (7,6,4,4,3,3,1)$,
the first $2$-Durfee rectangle of height~$s_2(\la) = 5$ and width~$3$, and the
second $2$-Durfee square of height~$t_2(\la) = 2$ and width~$0$.  Here the
remaining partitions are $\al = (4,3,1,1)$, $\be = (3,1)$,
and $\ga$ which is empty. In this case, we have $(2,2)$-rank
$r_{2,2}(\la) = \be_1 + \al_1 - \ga_1^\pr =  3+4-0 = 7$. }
\label{fig:durfee2}
\end{figure}

For $m \geq 0$, define an $m$\emph{-rectangle} to be a rectangle whose height
minus its width is $m$.
Define the {\em first $m$-Durfee rectangle} to be the largest $m$-rectangle
which fits in diagram~$[\la]$.  Denote by $s_m(\la)$ the
height of the first $m$-Durfee rectangle.  Define the {\em second $m$-Durfee rectangle}
to  be the largest $m$-rectangle which fits in diagram~$[\la]$ below
the first $m$-Durfee rectangle, and let $t_m(\la)$ be its height.
We will allow an $m$-Durfee rectangle to have width~0 but never height~0.
Finally, denote by~$\al$,~$\be$, and~$\ga$ the three partitions
to the right of, in the middle of and below the $m$-Durfee rectangles
(see Figures~\ref{fig:durfee} and~\ref{fig:durfee2}).  Notice that if $m > 0$ and we have an $m$-Durfee rectangle of width~$0$, as in Figure~\ref{fig:durfee2}, then $\ga$ must be the empty partition.

We define {\em $(2,m)$-rank}, $r_{2,m}(\la)$, of a partition~$\la$
by the formula:
$$
r_{2,m}(\la) \, := \, \be_1 + \al_{s_m(\la)-t_m(\la) - \be_1+1} - \ga_1^\pr \, .
$$
Note that $(2,0)$-rank is only defined for non-Rogers-Ramanujan
partitions because otherwise $\be_1$ does not exist, while
$(2,m)$-rank is defined for all partitions for all~$m >0$. Again, see Figures~\ref{fig:durfee} and~\ref{fig:durfee2} for examples.

Let~$\ch_{n,m,r}$ be the set of partitions of~$n$ with~$(2,m)$-rank~$r$.
In the notation above, $h(n,m,r) = \bigl|\ch_{n,m,r}\bigr|$.  Define~$\ch_{n,m, \leq r}$ and~$\ch_{n,m, \geq r}$ similarly.

\medskip


\subsection{Proof of the first symmetry}\label{sec-sym1} \
In order to prove
the first symmetry we present an involution~$\vp$ on~$\cp \sm \cq$
which preserves the size of partitions as well as their Durfee
squares, but changes the sign of the rank:
$$
\vp: \ch_{n,0,r} \to \ch_{n,0,-r} \, .
$$

Let~$\la$ be a partition with two Durfee square and partitions~$\al$,~$\be$, and~$\ga$ to the right of,
in the middle of, and below the Durfee squares.  This map~$\vp$ will preserve
the Durfee squares of~$\la$ whose sizes we denote by
$$
s = s_0(\la) \hspace{.66cm} \text{and} \hspace{.66cm} t = t_0(\la)\, .
$$

We will describe the action of~$\vp: \la \mapsto \wh \la$ by
first mapping~$(\al, \be, \ga)$
to a 5-tuple
of partitions~$(\mu,\nu,\pi,\rho,\si)$, and subsequently mapping that 5-tuple to
different triple~$(\wh{\al}, \wh{\be}, \wh{\ga})$ which goes to the right of, in the
middle of, and below the Durfee squares in~$\wh{\la}$.

\smallskip

\noindent
\begin{enumerate}
\item
First, let~$\mu = \be$.

\noindent Second, remove the following parts from~$\al$:~$\al_{s - t - \be_j +j}$ for~$1 \leq  j \leq t$.
Let~$\nu$ be the partition
comprising of parts removed from~$\al$ and~$\pi$ be the partitions comprising of the parts
which were not removed.

\noindent Third,
for~$1 \leq j \leq t$, let
$$
k_j = \mathrm{max}\{k \leq s - t \mid
\ga_j^\pr -k \geq \pi_{s -t -k +1}\} \, .
$$
Let~$\rho$ be the partition with
parts~$\rho_j = k_j$ and $\si$ be the partition with parts $\si_j = \ga_j^\pr - k_j$.
\item
First, let~$\wh{\ga}^\pr = \nu + \mu$ be the {\em sum} of partitions,
defined to have parts~$\wh\ga^\pr_j = \nu_j + \mu_j$.

\noindent Second, let~$\wh{\al} = \si \cup \pi$ be the {\em union} of
partitions, defined
as a union of parts in~$\si$ and~$\pi$.\footnote{Alternatively,
the union can be defined via the sum: $\si \cup \pi = (\si'+\pi')'$.}

\noindent Third, let~$\wh{\be} = \rho$.

\end{enumerate}

\smallskip

\begin{figure}[hbt]
\begin{center}
\psfrag{s}{$s$}
\psfrag{phi}{$\vp$}
\psfrag{la}{$\la$}
\psfrag{lp}{$\la^\pr$}
\psfrag{al}{$\al$}
\psfrag{be}{$\be$}
\psfrag{ga}{$\ga$}
\psfrag{gp}{$\ga'$}
\psfrag{pi}{$\pi$}
\psfrag{si}{$\si$}
\psfrag{ro}{$\rho$}
\psfrag{nu}{$\nu$}
\psfrag{mu}{$\mu$}
\psfrag{xx}{$\wh{\ga}$}
\psfrag{xp}{$\wh{\ga}^{\prime}$}
\psfrag{yy}{$\wh \al$}
\psfrag{zz}{$\wh \be$}
\psfrag{=}{$=$}
\psfrag{+'}{$+$}
\psfrag{'+}{$+$}
\psfrag{C}{$\cup$}
\epsfig{file=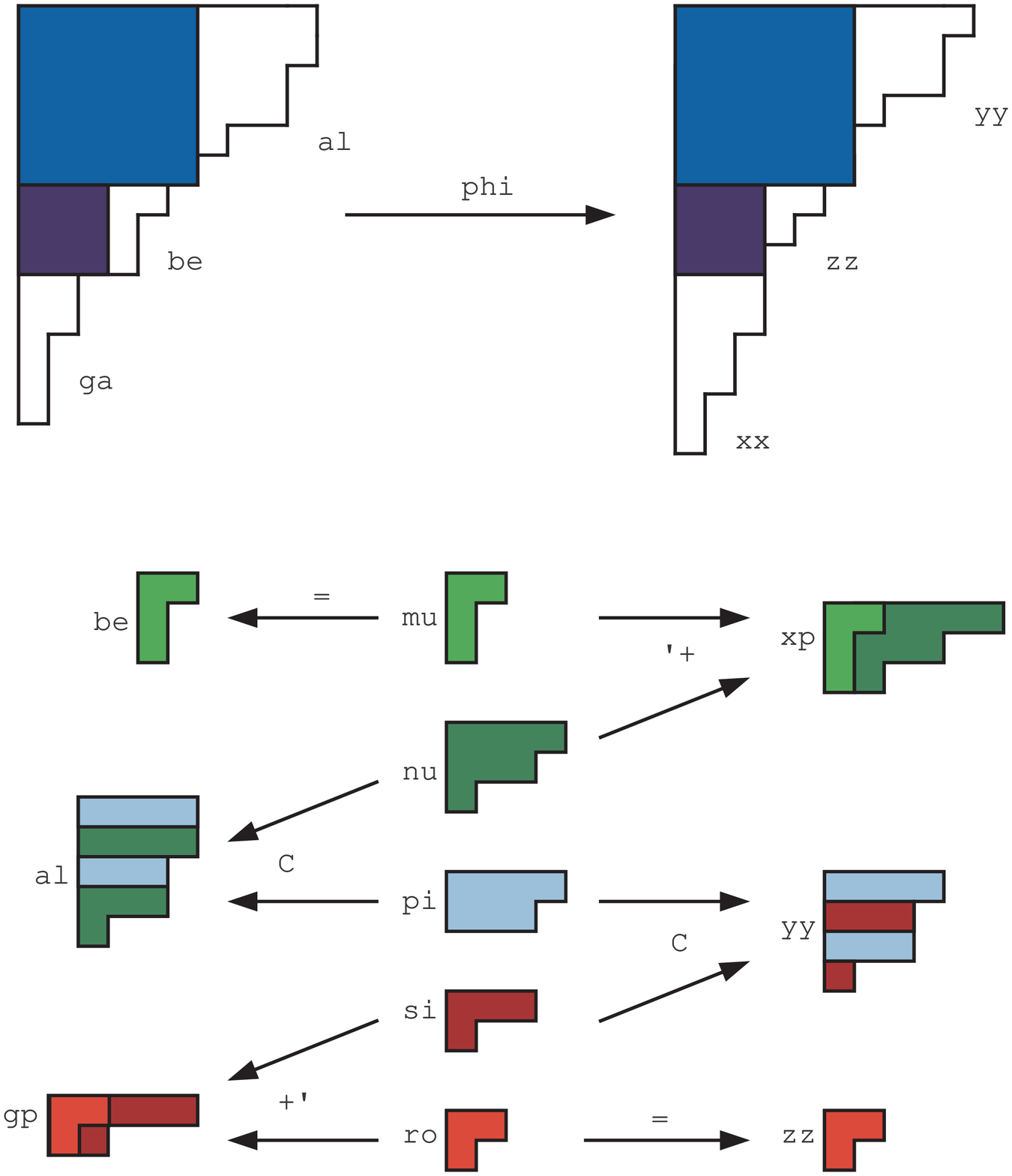,height=14.8cm}
\end{center}
\caption{An example of the first symmetry involution
$\vp: \la \mapsto \wh \la$, where $\la \in \ch_{n,0,r}$ and
$\wh \la \in  \ch_{n,0,-r}$ for $n= 71$, and $r=1$.
The maps are defined by the following rules:
$\be = \mu$, $\al = \nu \cup \pi$, $\ga^\pr = \si + \rho$,
while $\wh{\be} = \rho$, $\wh{\al} = \pi \cup \si$,
 $\wh{\ga}^\pr = \mu + \nu$.  Also,
 $\la = (10,10,9,9,7,6,5,4,4,2,2,1,1,1)$ and
$\wh \la = (10,9,9,7,6,6,5,4,3,3,3,2,2,1,1)$.
 }
\label{fig:first}
\end{figure}

\noindent
Figure~\ref{fig:first} shows an example of $\vp$ and the relation between
these two steps.

\begin{Remark}\label{rem-sym1}
{\rm
The key to understanding the map $\vp$ is the definition of $k_j$.
By considering $k=0$, we see that $k_j$ is defined for all $ 1 \leq j \leq t$.
Moreover, one can check that $k_j$ is the \emph{unique} integer~$k$
which satisfies
$$
(\dag) \qquad \pi_{s - t -k +1} \leq \ga_j^\pr -k \leq \pi_{s-t -k} \, .
$$
(We do not consider the upper bound for $k = s -t$.)
This characterization of $k_j$ can also be taken as its definition.
Equation~$(\dag)$ is used repeatedly in our proof of the next lemma.
}
\end{Remark}

\begin{lemma}
The map $\vp$ defined above is an involution.
\end{lemma}

\begin{proof}
Our proof is divided into five parts; we prove that

\smallskip

\textbf{(1)} \ $\rho$ is a partition, \qquad
\textbf{(2)} \ $\si$ is a partition, \qquad
\textbf{(3)} \ $\wh{\la} = \vp(\la)$ is a partition,

\smallskip

\textbf{(4)} \ $\vp^2$ is the identity map, \qquad
and \qquad \textbf{(5)} \
$r_{2,0}(\wh{\la}) = - r_{2,0}(\la)$.

\smallskip

\noindent \textbf{(1)} \
Considering the bounds ($\dag$) for $j$ and $j+1$,
we note that, if $k_j \leq k_{j+1}$, then
$$
\pi_{s -t -k_j +1} + k_j \, \leq \, \pi_{s -t -k_{j+1} +1} + k_{j+1} \, \leq \, \ga_{j+1}^\pr \leq \ga_j^\pr \, \leq \pi_{s -t -k_j} + k_j\, .
$$
This gives us
$$
\pi_{s -t -k_j +1} \leq
\ga_{j+1}^\pr -k_j \leq \pi_{s -t -k_j}
$$
and uniqueness therefore implies that $k_j = k_{j+1}$.
We conclude that $k_j \geq k_{j+1}$ and that $\rho$ is a partition.

\medskip

\noindent \textbf{(2)} \
If $k_j > k_{j+1}$, then we have $s -t -k_j +1  \leq  s -t -k_{j+1}$
and therefore
$$
\pi_{s-t -k_{j+1}} \leq \pi_{s - t -k_j +1} \, .
$$
Again, by considering ($\dag$) for $j$ and $j+1$, we conclude
that
$$
\ga_j^\pr -k_j \geq \ga_{j+1}^\pr -k_{j+1} \, .
$$
If $k_j = k_{j+1}$, then we simply need to recall that $\ga^\pr$ is a partition to see that
$$
\ga_j^\pr -k_j \geq \ga_{j+1}^\pr -k_{j+1} \, .
$$
This implies that~$\si$ is a partition.

\medskip

\noindent \textbf{(3)} \
By their definitions, it is clear that~$\mu$,~$\nu$, and~$\pi$ are partitions.
Since we just showed that~$\rho$ and~$\si$ are all partition, it follows that~$\wh{\al}$,~$\wh{\be}$,
and~$\wh{\ga}$ are also partitions.
Moreover, by their definitions,
we see that~$\mu$,~$\nu$, and~$\si$ have at most~$t$ parts,~$\pi$ has at most~$s-t$,
and~$\rho$ has at most $t$ parts each of which is less than or equal to~$s-t$.
This implies that~$\wh{\al}$ has at most~$s$ parts,~$\wh{\be}$ has at most~$t$ parts
each of which is less than or equal to~$s-t$, and~$\wh{\ga}^\pr$ has parts at most~$t$.
Therefore,~$\wh{\al}$,~$\wh{\be}$, and~$\wh{\ga}$ fit to the right of, in the middle of, and below
Durfee squares of sizes~$s$ and~$t$ and so~$\vp(\la)$ is a partition.

\medskip

\noindent \textbf{(4)} \
We will apply~$\vp$ twice to a non-Rogers-Ramanujan partition~$\la$
with~$\al$,~$\be$, and~$\ga$ to the right of, in
the middle of, and below its two Durfee squares.
As usual, let~$\mu$,~$\nu$,~$\pi$,~$\rho$,~$\si$ be the partitions occurring
in the intermediate stage of the first application of~$\vp$ to~$\la$ and
let~$\wh{\al}$,~$\wh{\be}$,~$\wh{\ga}$
be the partitions to the right of, in the middle of, and below
the Durfee squares of~$\wh{\la}=\vp(\la)$.
Similarly,
let~$\wh{\mu}$,~$\wh{\nu}$,~$\wh{\pi}$,~$\wh{\rho}$,~$\wh{\si}$ be the
partitions occurring in the intermediate stage of the second application
of~$\vp$ and let~$\al^\ast$,~$\be^\ast$, and~$\ga^\ast$ be the partitions
to the right of, in the middle of and below the Durfee squares
of~$\vp^2(\la) = \vp(\wh{\la})$.

We need several observations. \,
First, note that~$\wh{\mu} = \wh{\be} = \rho$.  \, Second, by ($\dag$)
we have:
$$
\pi_{s - t -k_j +1} \leq \ga_j^\pr -k_j = \si_j \leq \pi_{s-t -k_j} \,.
$$
Since~$\si$ is a partition, this implies that~$\wh{\al}_{s - t -k_j +j} = \si_j$.
On the other hand, since~$\wh{\be}_j = \rho_j = k_j$, the map~$\vp$ removes the
rows~$\wh{\al}_{s - t -k_j +j} = \si_j$ from~$\wh{\al}$.
From here we conclude that~$\wh{\nu} = \si$ and~$\wh{\pi} = \pi$.
\, Third, define
$$
\wh{k}_j = \mathrm{max}\{\wh{k} \leq s - t \mid
\ga_j^\pr -\wh{k} \geq \pi_{s -t -\wh{k} +1}\} \, .
$$
By Remark~\ref{rem-sym1}, we know that~$\wh{k}_j$ as above is the
unique integer~$\wh{k}$ which satisfies:
$$
\wh{\pi}_{s - t -\wh{k} +1} \leq \wh{\ga}_j^\pr -\wh{k} \leq \wh{\pi}_{s-t -\wh{k}} \, .
$$
On the other hand, recall that~$\wh{\ga}_j^\pr = \mu_j + \nu_j$ and~$\be_j = \mu_j$.
This implies~$\wh{\ga}_j^\pr - \be_j = \nu_j$.
Also, by the definition of~$\nu$, we have $\nu_j = \al_{s-t-\be_j+j}$.
Therefore, by the definition of~$\pi$, we have:
$$
\pi_{s - t - \be_j +1} \leq \al_{s - t - \be_j +j} = \nu_j = \wh{\ga}_j^\pr - \be_j \leq
\pi_{s-t -\be_j} \, .$$
Since,~$\wh{\pi} = \pi$, by the uniqueness in Remark~\ref{rem-sym1}
we have~$\wh{k}_j = \be_j = \mu_j$. This implies
that~$\wh{\rho} = \mu$ and~$\wh{\si} = \nu$.

Finally, the second step of our bijection gives $\al^\ast = \nu \cup \pi = \al$,
$\be^\ast = \mu = \be$, and~$(\ga^\ast)^\pr = \rho + \si = \ga^\pr$.
This implies that~$\vp^2$ is the identity map.

\medskip

\noindent \textbf{(5)} \
Using the results from \textbf{(4)}, we have:
$$r_{2,0}(\la) = \be_1 + \al_{s-t - \be_1+1} - \ga_1^\pr
= \mu_1 + \nu_1 - \rho_1 - \si_1\,.$$
On the other hand,
$$r_{2,0}(\wh{\la}) = \wh{\be}_1 + \wh{\al}_{s-t -\wh{\be}_1+1} -
\wh{\ga}_1^\pr = \rho_1 + \si_1 - \mu_1  - \nu_1 \, .$$
We conclude that $r_{2,0}(\wh{\la}) = - r_{2,0}(\la)$.
\end{proof}

\medskip


\subsection{Proof of the second symmetry}\label{sec-sym2} \
In order to prove the second symmetry we present a bijection
$$
\psi_{m,r}: \ch_{n,m,\leq -r} \to \ch_{n-r-2m-2,m+2,\geq -r} \, .
$$
This map will only be defined for~$m,r > 0$ and for $m = 0$ and $r \geq 0$ and in both of these cases
the first and second~$m$-Durfee rectangles of a partition $\la \in \ch_{n,m,\leq -r}$
have non-zero width.  For~$m =0$, $(2,0)$-rank is only defined for partitions in~$\cp \sm \cq$ which by definition have two Durfee squares of non-zero width.  For $m>0$, since we also have~$r>0$, a partition $\la \in \ch_{n,m,\leq -r}$ must have
$$
r_{2,m}(\la) = \be_1 + \al_{s_m(\la) - t_m(\la) - \be_1+1} -\ga_1^\pr \leq -r < 0 \, .
$$
This forces $\ga_1^\pr > 0$ and so both $m$-Durfee rectangles must have non-zero width.

We describe the action of~$\psi := \psi_{m,r}$ by giving
the sizes of the Durfee
rectangles of~$\wh{\la} := \psi_{m,r}(\la) = \psi(\la)$
and the partitions~$\wh{\al}$,~$\wh{\be}$, and~$\wh{\ga}$
which go to the right of, in the middle of,
and below those Durfee rectangles in~$\wh{\la}$.

\smallskip

\noindent
\begin{enumerate}
\item
If~$\la$ has two $m$-Durfee rectangles of height
$$
s := s_m(\la) \hspace{.66cm} \text{and} \hspace{.66cm} t := t_m(\la)
$$
then $\wh{\la}$ has two $(m+2)$-Durfee rectangles of height
$$
s^\pr:= s_{m+2}(\wh{\la}) = s+1 \hspace{.66cm} \text{and} \hspace{.66cm}
t^\pr:= t_{m+2}(\wh{\la}) = t+1 \, .
$$
\item
Let
$$
k_1 = \mathrm{max}\{k \leq s -t \mid \ga_1^\pr -r -k \geq \al_{s-t -k+1}\} \, .
$$
Obtain $\wh{\al}$ from $\al$ by adding a new part of size $\ga_1^\pr -r -k_1$,
$\wh{\be}$ from $\be$ by adding a new part of size $k_1$,
and $\wh{\ga}$ from $\ga$ by removing its first column.
\end{enumerate}

\smallskip

\begin{figure}[hbt]
\begin{center}
\psfrag{a}{$\ga_1'$}
\psfrag{k1}{$k_1$}
\psfrag{gk}{$\ga_1'-r-k_1$}
\psfrag{la}{$\la$}
\psfrag{hla}{$\wh \la$}
\psfrag{s}{$s$}
\psfrag{t}{$t$}
\psfrag{s'}{$s'$}
\psfrag{s''}{$s''$}
\psfrag{t'}{$t'$}
\psfrag{t''}{$t''$}
\psfrag{map}{$\psi_{m,r}$}
\epsfig{file=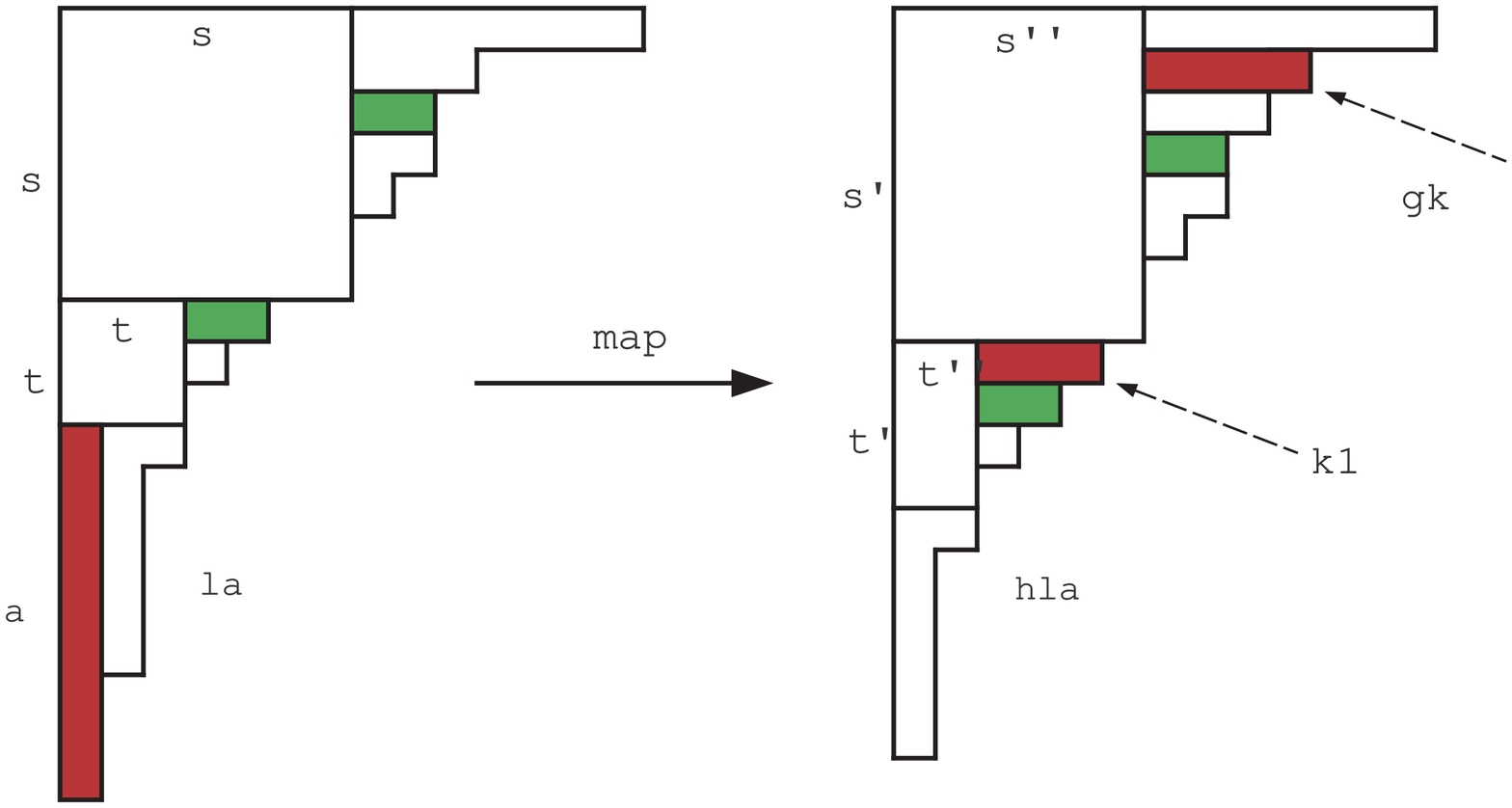,height=7.2cm}
\end{center}
\caption{An example of the second symmetry bijection
$\psi_{m,r}: \la \mapsto \wh \la$, where
$\la \in \ch_{n,m,\le -r}$,
$\wh \la \in  \ch_{n',m+2,\ge -r}$, for $m=0$, $r=2$,
$n=92$, and $n' = n -r-2m-2 = 88$.
Here  $r_{2,0}(\la) = 2+2-9 = -5 \le -2$ and
$r_{2,2}(\wh \la) = 3+4-6 = 1 \ge -2$, where
$\la = (14,10,9,9,8,7,7,5,4,3,3,2,2,2,2,2,1,1,1)$
 and
$\wh \la = (13,10,9,8,8,7,6,6,5,4,3,2,2,1,1,1,1,1)$.
Also, $s=7$, $s' = s+1 = 8$, $s'' =s' -m -2= 6$,
$t = 3$, $t'=4$, $t'' = 2$,
$\ga_1' = 9$, $k_1 = 3$, and $\ga_1' -r-k_1 = 4$.
}
\label{fig:second}
\end{figure}

\noindent
Figure~\ref{fig:second} shows an example of the
bijection~$\psi = \psi_{m,r}$.

\begin{Remark}\label{rem-sym2}
{\rm
As in Remark~\ref{rem-sym1},
by considering~$k = \be_1$ we see that~$k_1$ is defined and
indeed we have~$k_1 \geq \be_1$.
Moreover, it follows from its definition that~$k_1$ is the
{\em unique}~$k$ such that
$$ (\ddag) \qquad
\al_{s - t -k +1} \leq \ga_1^\pr -r-k \leq
\al_{s-t -k} \, .
$$
(If $k = s - t$ we do not consider the upper bound.)
}
\end{Remark}

\begin{lemma}
The map $\psi = \psi_{m,r}$ defined above is a bijection.
\end{lemma}

\begin{proof}
Our proof has four parts:

\smallskip
\textbf{(1)} \ we prove that~$\wh{\la} = \psi(\la)$ is a partition,

\textbf{(2)} \ we prove that the size of~$\wh{\la}$ is~$n-r-2m-2$,

\textbf{(3)} \ we prove that $r_{2,m+2}(\wh{\la}) \geq -r$, and

\textbf{(4)} \ we present the inverse map~$\psi^{-1}$.

\medskip

\noindent  \textbf{(1)} \
To see that~$\wh{\la}$ is a partition we simply have to note that since $\la$ has $m$-Durfee rectangles of non-zero width, $\wh{\la}$ may have $(m+2)$-Durfee rectangles of width $s-1$ and $t-1$.  Also, the partitions $\wh{\al}$ and $\wh{\be}$ have at most $s+1$ and $t+1$ parts, respectively, while the partitions $\wh{\be}$ and $\wh{\ga}$ have parts of size at most $s-t$ and $t-1$, respectively.  This means that they can sit to the right of, in the middle of, and below the two $(m+2)$-Durfee rectangles of $\wh{\la}$.

\medskip

\noindent \textbf{(2)} \
To prove that the above construction gives a partition~$\wh{\la}$
of~$n-r-2m-2$, note that the sum of the sizes of the rows added to~$\al$
and~$\be$ is~$r$ less
than the size of the column removed from~$\ga$, and that
both the first and second~$(m+2)$-Durfee rectangles of~$\wh{\la}$
have size~$m+1$ less
than the size of the corresponding~$m$-Durfee rectangle of~$\la$.

\medskip

\noindent  \textbf{(3)} \
By Remark~\ref{rem-sym2}, the part we inserted into~$\be$ will
be the largest part of the resulting
partition, i.e.~$\wh{\be}_1 = k_1$.
By equation~($\ddag$) we have:
$$
\al_{s - t -k_1 +1} \leq \ga_1^\pr -r-k_1 \leq \al_{s-t-k_1} \, .
$$
Therefore, we must have:
$$
\wh{\al}_{s^\pr-t^\pr -\wh{\be}_1 +1} =
\wh{\al}_{s-t -k_1 +1} = \ga_1^\pr -r-k_1 \, .
$$
Indeed, we have chosen~$k_1$ in the \emph{unique} way so that the rows we insert into~$\al$ and~$\be$ are~$\wh{\al}_{s^\pr-t^\pr -\wh{\be}_1 +1}$
and~$\wh{\be}_1$ respectively.

Having determined~$\wh{\al}_{s^\pr-t^\pr -\wh{\be}_1 +1}$
and~$\wh{\be}_1$ allows us to bound the~$(2,m+2)$-rank of~$\wh{\la}$ :
$$
r_{2,m+2}(\wh{\la}) = \wh{\al}_{s^\pr-t^\pr -\wh{\be}_1 +1} + \wh{\be}_1 - \wh{\ga}_1^\pr \, = \,
\ga_1^\pr -r -k_1 +k_1 - \wh{\ga}_1^\pr \, \geq \, -r\,,
$$
where the last inequality follows since~$\wh{\ga}_1^\pr$ is the size of the second column of~$\ga$ whereas~$\ga_1^\pr$ is the size of the first column of~$\ga$.

\medskip

\noindent  \textbf{(4)} \
The above characterization of~$k_1$ also shows us that to
recover~$\al$,~$\be$, and~$\ga$ from~$\wh{\al}$,~$\wh{\be}$ and~$\wh{\ga}$,
we remove part~$\wh{\al}_{s^\pr-t^\pr -\wh{\be}_1 +1}$ from~$\wh{\al}$,
remove part~$\wh{\be}_1$ from~$\wh{\be}$, and add a column of
height~$\, \wh{\al}_{s^\pr-t^\pr -\wh{\be}_1 +1} + \wh{\be}_1 +r \,$
to~$\wh{\ga}$.
Since we can also easily recover the sizes of the previous $m$-Durfee rectangles,
we conclude that~$\psi$ is a bijection between the desired sets.
\end{proof}

\bigskip

\section{Final remarks} \label{sec:final}

\subsection{}\label{sec:fin-1}  Of the many proofs of
Rogers-Ramanujan identities only a few can be honestly
called ``combinatorial''.   We would like to single
out~\cite{A3} as an interesting example.  Perhaps, the most
important combinatorial proof was given by Schur in~\cite{S}
where he deduced his identity by a direct involutive argument.
The celebrated bijection of Garsia and Milne~\cite{GM} is based on this
proof and the involution principle.
In~\cite{BZ1}, a different involution principle proof
was obtained (see also~\cite{BZ3}) based on a short
proof of Bressoud~\cite{B}.  We refer to~\cite{P1} for
further references, historical information, and combinatorial proofs of other partition identities.

\subsection{}\label{sec:fin-1.5}Dyson's rank $r_1(\la) = \la_1 - \la_1^\pr$ was
defined in~\cite{D1} for the purposes of finding a combinatorial
interpretation of Ramanujan's congruences.  Dyson used
the rank to obtain a simple combinatorial proof of
Euler's pentagonal theorem in~\cite{D2}
(see also~\cite{D3,P2}).  It was shown in~\cite{P2}
that this proof can be converted into a direct involutive
proof, and such a proof in fact coincides with the involution
obtained by Bressoud and Zeilberger~\cite{BZ2}.

Roughly speaking, our proof of Schur's identity is a
Dyson-style proof with a modified Dyson's rank, where the
definition of the latter was inspired by~\cite{BZ1, BZ2, BZ3}.
Unfortunately, reverse engineering the proofs in~\cite{BZ3}
is not straightforward due to the complexity of that paper.
Therefore, rather than giving a formal connection, we will only
say that, for some~$m$ and~$r$, our map~$\psi_{m,r}$ is similar to the
maps~$\vp$ in~\cite{BZ1} and~$\it{\Phi}$ in~\cite{BZ3}.

It would be interesting to extend our Dyson-style proof to
the generalization of Schur's identity found in~\cite{GIS}.
This would give a new combinatorial proof of the generalizations
of the Rogers-Ramanujan identities found in that paper and, in a
special case, provide a new combinatorial proof
of the second  Rogers-Ramanujan identity (see
e.g.~\cite{A1,A4,H,P1}).

\subsection{}\label{sec:fin-2}  The idea of using
iterated Durfee squares to study the Rogers-Ramanujan identities
and their generalizations is due to Andrews~\cite{A2}.
The $(2,m)$-rank of a partition
is a special case of a general (but more involved)
notion of $(k,m)$-rank which is presented in~\cite{Bo}.
It leads to combinatorial proofs of some of
Andrews' generalizations of Rogers-Ramanujan identities mentioned above.

Garvan~\cite{G} defined a generalized notion of a rank to partitions with iterated Durfee squares,
that is different from ours, but still satisfies
equation~$(\maltese)$ (for $m=0$).
In~\cite{BG}, Berkovich and Garvan asked for a Dyson-style proof of~$(\maltese)$ but
unfortunately, they
were unable to carry out their program in
full as the combinatorial symmetry they obtain
seem to be hard to establish bijectively.  (This
symmetry is somewhat different from our second symmetry.)
The first author was able to
relate the two generalizations of rank
by a
bijective argument.  This also appears in~\cite{Bo}.

\subsection{}\label{sec:fin-2.5} Yet another generalization of Dyson's
rank was kindly brought to our attention by George Andrews. The notion of
successive rank can also be used to give a combinatorial proof of the
Rogers-Ramanujan identities and their generalizations by a sieve argument
(see~\cite{A5,B1}).
However, this proof involves a different combinatorial
description of the partitions on the left hand side of the
Rogers-Ramanujan identities than the proof presented here.

\subsection{}\label{sec:fin-3}  Finally, let us note that the
Jacobi triple product identity has a combinatorial proof
due to Sylvester (see~\cite{P1,W}).
We refer to~\cite{A0} for an elementary algebraic proof.

Also, while our proof is mostly combinatorial it is by no
means a direct bijection.  The quest for a direct bijective
proof is still under way, and as recently as this year
Zeilberger lamented on the lack of such proof~\cite{Z}.
The results in~\cite{P3} seem to discourage any future
work in this direction.

\vskip.6cm
\noindent
{\bf Acknowledgments.} \ The authors are grateful to
George Andrews and Richard Stanley for their support
and encouragement of our studies of partition identities.
The first author was supported by NSERC(Canada) and the second author by the NSA and the NSF.


\end{document}